\documentclass[12pt]{amsart}

\usepackage{amsmath, amsthm, amsfonts}
\newtheorem{proposition}{Proposition}
\newtheorem{theorem}{Theorem}

\newcommand{\ze}{\mathbb{Z}}
\newcommand{\re}{\mathbb{R}}

\author{Gregorio Moreno}

\address{\noindent Universit\'e Paris 7, Math\'ematiques,
 case 7012, 2, Place Jussieu,
75251 Paris, France
{\rm and}
\noindent Facultad de Matem\'aticas, Pontificia Universidad Cat\'olica de Chile,
Vicu\~na Mackena 4860, Macul, Chile
}

\email{gregorio.random@gmail.com}

\thanks{Partially supported by Beca Conicyt-Ambassade de France and CNRS, UMR $7599$ "Probabilit\'es et
Mod\`eles Al\'eatoires".}

\date{}

\title[Environment Seen by the Particle
for Directed Polymers]{Convergence of the law of the Environment Seen by the Particle
for Directed Polymers in Random Media in the $L^2$ region}

\begin{document}

\begin{abstract} We consider the model of Directed Polymers
in an i.i.d.  gaussian or bounded Environment \cite{IS, CH, CSY} in the $L^2$ region. We prove the convergence 
of the law of the environment seen by the particle. 

As a main technical step, we establish
a lower tail concentration inequality for the partition function for bounded environments. Our proof is based
on arguments developed by Talagrand in the context of the Hopfield Model \cite{T}. This improves in some sense
a concentration inequality obtained by Carmona and Hu for gaussian environments \cite{CH}. We use this and 
a Local Limit Theorem \cite{S, V} to prove the $L^1$ convergence of the density of the law
of the environment seen by the particle with respect to the product measure.

\end{abstract}

\maketitle

\section{Introduction}

We consider the following model of directed polymers in a random media: let $d\geq 3$ and let $\{ \eta(t,x):\, t\in \ze,\, x\in \ze^d\}$ denotes a family of real variables. We call it the environment.
For $y\in \ze^d$, let $P_y$ be the law of the simple symmetric nearest neighbor random walk on $\ze^d$ starting at $y$.
For $L\in \mathbb{N}$, we denote the space of nearest neighbor paths of length $L$ by $\Omega_L$, i.e.

\begin{eqnarray*}
\Omega_L := \left\lbrace s: \{0,1,2,...,L\} \to \ze^d,\, |s_{t+1} - s_t| = 1,\, \forall \, t=0,...,L-1 \right\rbrace,
\end{eqnarray*}

\noindent where $|\cdot|$ is the euclidean norm. Denote by $\omega$ the canonical process on
$\Omega_{L}$, i.e., $\omega_t(s) = s_t$ for $s\in \Omega_L$ and $t=0,1,...,L$. 

For a fixed configuration $\eta$, $y\in \ze^d$, $M,\, N \in \ze,\, M < N$ and $0<\beta<+\infty$, we can define the (quenched) law of the polymer in environment $\eta$, inverse temperature $\beta$, based on $(M,y)$ and time horizon $N$ (or simply the polymer measure): for all $s \in \Omega_{N-M}$, 

\begin{eqnarray}
 \mu^{y}_{M,N}(\omega = s)= \frac{1}{Z^y_{M,N}} \exp \left\{\beta\!\! \sum^{N-M}_{t=1}\!\!\!\eta(t+M,s_{t})\right\} P^y(\omega = s),
\end{eqnarray}

\vspace{2ex}
\noindent where

\begin{eqnarray}
 Z^y_{M,N} = P^y \left[%
 \exp \left\{ \beta \!\!\! \sum^N_{t=M+1}\!\!\eta(t+M,\omega_{t})\right\}%
\right],
\end{eqnarray}

\noindent is usually called the partition function. $\mu^y_{M,N}$ can be though as a measure on directed
paths, starting from $y$ at time $M$ and ending at time $N$. When $M=0$ and $y=0$, we will sometimes not
write them.

We endow the space $E=\re^{\ze^{d+1}}$ with the product $\sigma$-algebra and a product measure $Q$.
We can then see $\eta$ as a random variable with values in $E$. In the following we will deal with
laws $Q$ which marginals are gaussian or have bounded support. We are interested in the properties
of the polymer measure for a typical configuration $\eta$.

A quantity of special interest is $\lambda(\beta)= \log Q(e^{\beta \eta(0,0)})$. We can compute it
for simple laws of the environment. For example, in the case of a standard gaussian law, we have
$\lambda(\beta)=\beta^2/2$.

The model of Directed Polymers in Random Environment as been extensively studied over the
last thirty years \cite{IS, CH, CSY, CV, V2}. It turns out that the behavior of the polymer depends strongly on the
temperature. Here we will be concerned by a region of very high temperature, or equivalently, we will focus
on very small values of $\beta$. Let state this more precisely: we can define the normalized partition function,
\begin{eqnarray}
 W_N = Z_N \exp \{ -N\lambda(\beta)\},
\end{eqnarray}

\noindent  It follows easily that this is
a positive martingale with respect to the filtration $\mathcal{H}_N=\sigma(\eta(t,x): t\leq N)$. Then, by 
classical arguments, we can prove that it converges when $N \to +\infty$ to a non-negative random variable
$W_{+\infty}$ that satisfies the following zero-one law:
\begin{eqnarray*}
 Q(W_{+\infty}>0) = 0\,\, {\rm or}\,\, 1.
\end{eqnarray*}

\noindent In the first situation, we say that strong disorder holds. In the second, we say that
weak disorder holds. The behavior is qualitatively different in each situation \cite{CSY}. Let recall
an argument showing that the weak disorder region is nontrivial, at least for $d\geq 3$ and
a wide class of environments. Observe first that in order to prove weak disorder, it is enough to
obtain some uniform integrability condition on the martingale $(W_N)_N$. Indeed, uniform integrability
implies that $Q(W_{+\infty})=1$, and then the zero-one law gives us that $Q(W_{+\infty}>0)=1$.
The easiest situation occurs when the martingale is bounded in $L^2$: Define

\begin{equation}\nonumber
 \pi_d = P_0^{\otimes 2}\left[\omega_t = \tilde{\omega}_t \, \text{ for some $t$}\right].
\end{equation}

\noindent We see that $\pi_d<1$ only for $d\geq 3$. Let's perform the following elementary calculation:

\begin{eqnarray}\label{L2}
 Q(W^2_{N})&=&Q P^{\otimes 2}\left[\exp \left\{ \beta \sum^N_{t=1} \eta(t,\omega_t)+
                                    \beta \sum^N_{t=1} \eta(t,\tilde{\omega}_t)
                                    - 2N \lambda(\beta) \right\} \right]\\ \nonumber
&=&e^{-2N \lambda} P^{\otimes 2}\left[ \prod^N_{t=1}Q(e^{ 2\beta \eta(t,\omega_t){\bf 1}_{\omega_t=\tilde{\omega}_t}}) \times Q(e^{\beta \eta(t,\tilde{\omega}_t){\bf 1}_{\omega_t \neq \tilde{\omega}_t}})^2\right]\\ \nonumber
           &=&P^{\otimes 2}\left[e^{ (\lambda(2\beta)-2\lambda(\beta))L_N(\omega,\tilde{\omega})}\right],
\end{eqnarray}

\noindent where $L_N(\omega,\tilde{\omega})= \sum^N_{t=1}{\bf 1}(\omega_t = \tilde{\omega}_t)$. We observe that
$L_{+\infty}(\omega,\tilde{\omega})$ has a geometric distribution with parameter $\pi_d$, so that $Q(W^2_N)$ is
uniformly bounded in $N$ for 

\begin{eqnarray}\label{condL2}
 \lambda(2\beta)-2\lambda(\beta)< \log (1/\pi_d),
\end{eqnarray}

\noindent and uniform integrability follows. We call this the $L^2$ condition. In the general case, we can prove that, in fact, such a condition
holds for sufficiently small values of $\beta$ (we can see this directly for gaussian environments). The range of
values for which (\ref{condL2}) holds is called the $L^2$ region. 

We know \cite{IS, B, AZ, SZ} that in the $L^2$ region, the polymer is diffusive. Indeed, we have the following 
(quenched) invariance principle:

\begin{theorem}\label{invariance}\cite{AZ, B, IS, SZ}
 In the $L^2$ region, for $Q$-almost every environment, we have that,
under the polymer measure,
\begin{eqnarray*}
 \frac{1}{\sqrt{N}}\, \omega(Nt)
\end{eqnarray*}

\noindent converges in law to a brownian motion with covariance matrix $1/d\, I$, where
$I$ is the identity matrix in dimension $d$.
\end{theorem}
\smallskip

\noindent More generally, in the full weak disorder region, a slightly weaker result holds (\cite{CY}).
The situation is quite different and more subtle in the strong disorder region. In that case,
large values of the medium attract the path of the polymer, so that a localization phenomenon arises
(for more information, see \cite{CH, CY}, and also \cite{V2} for milder assumptions).

\vspace{3ex}

In this article, we are concerned with another kind of (still related) result, namely the convergence of
the law of the environment viewed by the particle. The environment viewed by the particle is a
process $\eta_N$ with values in $E$, defined by

\begin{eqnarray}\label{esbp}
 \eta_N(t,x)= \eta(N+t,\omega_N+x), \quad t\in \ze,\, x\in \ze^d.
\end{eqnarray}

\noindent where $\omega$ follows the law $\mu^0_{0,N}$. Let's denote by $Q_N$ the law of this process at time $N$. 
We can see that

\begin{eqnarray}\label{density}
 \frac{dQ_N}{dQ}=\sum_{x}\mu^x_{-N,0}(\omega_N=0).
\end{eqnarray}

\noindent Indeed, take a bounded measurable function $f:E \to \re$. Then, by
translation invariance,

\begin{eqnarray}\nonumber
 Q \left(\mu^0_{0,N}f(\eta_N)\right)&=& Q \left( \sum_x  f(\eta(N+\cdot, x+\cdot))\mu^0_{0,N}(\omega_N=x)\right)\\
 \nonumber
&=&Q \left( f(\eta) \sum_x  \mu^x_{-N,0}(\omega_N=0)\right)\\
\label{calcdensity}
&=& Q_N\left( f(\omega)\right).
\end{eqnarray}

\noindent For a wide variety of models, the convergence of $Q_N$ has been a powerful tool for proving
invariance theorems. Obviously, it is not the case here because we have already the invariance
principle at hand. However, we can see this result as completing the picture in the $L^2$ region.
In the best of our knowledge, the point of view of the particle has never been studied before in
the literature for directed polymers. We can state the principal result of this article:
let $\overleftarrow{\eta}(t,x):= \eta(-t,x)$, $\overleftarrow{W}_N(\eta):= W_N(\overleftarrow{\eta})$
and let $\overleftarrow{W}_{+\infty}$ be the (almost sure) limit of $\overleftarrow{W}_N$ when
$N$ tends to infinity.

\vspace{3ex}




\begin{theorem}\label{convergence}
 In the $L^2$ region,

\begin{equation}\label{q_N}
 q_N:=\frac{dQ_N}{dQ} \longrightarrow \overleftarrow{W}_{+\infty} \times e^{\beta \eta(0,0)-\lambda(\beta)},\quad
  {\rm as}\,\, {N\to +\infty},
\end{equation}

\noindent where the limit is in the $L^1(Q)$ sense. 


\noindent In other words, $Q_N$ converges in the total variation distance to a probability measure
$Q_{+\infty}$ such that

\begin{eqnarray*}
 \frac{dQ_{+\infty}}{dQ}= \overleftarrow{W}_{+\infty} \times e^{\beta \eta(0,0)-\lambda(\beta)}.
\end{eqnarray*}

\end{theorem}
\smallskip




\noindent Much in the same spirit, we consider the law of the environment seen by the 
particle at an intermediate time $N$ under $\mu^0_{0,N+M}(\cdot)$. Formally, this new
environment is defined as the field $\eta_{N,M} \in E$ with

\begin{eqnarray*}
\eta_{N,M}(t,x)=\eta (N+t, \omega_N+x),\quad t\in \ze,\, x\in \ze^d,
\end{eqnarray*} 

\noindent where $\omega$ is taken from $\mu^0_{0,N+M}$. Following the argument of (\ref{calcdensity}), 
its density with respect to $Q$ is easily
seen to be

\begin{eqnarray}\label{densityNM}
q_{N,M} = \sum_{x}\mu^x_{-N,M}(\omega_{N}=0).
\end{eqnarray}

\begin{theorem}\label{convergence2}
In the $L^2$ region, we have

\begin{eqnarray*}
q_{N,M} \longrightarrow \overleftarrow{W}_{+\infty} \times
e^{\beta \eta(0,0)-\lambda(\beta)} \times W_{+\infty},\quad {\rm as}\,\, M,\,N\to +\infty,
\end{eqnarray*}

\noindent where the limit is in the $L^1(Q)$ sense.
\end{theorem}
\smallskip

\noindent A statement similar to Theorem $2$ may be found in Bolthausen and Sznitman \cite{BS} for directed random walks in
random environment. In their context, $q_N$ is a martingale. In our case, the denominator in the definition of
the polymer measure is quite uncomfortable as it depends on the whole past, so no martingale property should be 
expected for $q_N$. Let us mention that, we can cancel this denominator by multiplying each term of the sum in
(\ref{q_N}) by $W^x_{-N,0}$. This defines a martingale sequence that converges (almost surely!) to the same limit as $q_N$. 
In other terms, let $\mathcal{G}_N= \sigma(\eta(t,x):\, -N \leq t)$. Then, there exists a unique law $Q_{+\infty}$ on $E$ such that for any $A \in \mathcal{G}_N$,

\begin{eqnarray*}
Q_{+\infty}(A) = Q({\bf 1}_A \overleftarrow{W}_N e^{\beta \eta(0,0)-\lambda(\beta)}).
\end{eqnarray*}

\noindent The density of this law coincides with the limit in Theorem $2$ but we emphasize
that our strategy of proof here is completely different.

\vspace{2ex}

Let us mention that, in the weak disorder region, an infinite time horizon polymer measure
has been introduced in \cite{CY}. For each realization of the environment, it defines a Markov 
process with (inhomogeneous) transition probabilities given by

\begin{eqnarray*}
  \mu^{\beta}_{\infty}(\omega_{N+1}=x+e|\omega_N=x)=\frac{1}{2d} \, e^{\beta \eta(N+1,x+e)-\lambda(\beta)}\frac{W_{+\infty}(N+1,x+e)}{W_{+\infty}(N,x)},
\end{eqnarray*}

\noindent for all $x\in \ze^d$ and $|e|=1$, where for $k\in \ze$ and $y\in \ze^d$, 

\begin{eqnarray*}
W_{+\infty}(k,y)(\eta) = W_{+\infty}(\eta(k+\cdot,y+\cdot)).
\end{eqnarray*}

\noindent We can define the process of the environment seen by the particle for this
polymer measure by mean of the formula (\ref{esbp}), where, this time, $\omega_N$ 
is taken from $\mu^{\beta}_{+\infty}$. A simple modification of the proof of the Theorem $2$
suffices to show that the density of the environment seen by the particle converges,
as $N$ tends to $+\infty$, to the same limiting density as in Theorem $3$. The additional
term $W_{+\infty}$ arises from the very definition of the infinite time horizon measure.

\vspace{3ex}

We will prove the following lower tail concentration inequality. A similar statement as already
been proved by Carmona and Hu \cite{CH}, Theorem $1.5$, in the context of gaussian environments. As we will consider
bounded environments, there is no loss of generality if we assume that $\eta(t,x)\in [-1,1]$ for all
$t\in \ze,\, x\in \ze^d$.

\begin{proposition}\label{concentration}
 In the $L^2$ region, we can find $C>0$ such that,

\begin{equation}
 Q \left( \log Z_N \leq N\lambda - u \right) \leq C \exp \left\lbrace-\frac{u^2}{16 C \beta^2}\right\rbrace,
 \quad \forall \, N\geq 1,\, \forall \, u>0.
\end{equation}

\end{proposition}

\noindent It would be interesting to extend this result to a larger part of the
weak disorder region. This is indeed a major obstacle to extend our Theorems to 
larger values of $\beta$.

Another cornerstone in the proof of Theorem \ref{convergence} is a Local Limit Theorem
\cite{S, V} that we will describe later. Again, this result is available in the $L^2$ region only. 

\vspace{3ex}

We will now introduce some obvious notation that will be useful in what follows:
For $x\in \ze^d$ and $M<N \in \ze$, we write

\begin{eqnarray}
 W_{M,N}(x)=Z^x_{M,N}e^{-(N-M)\lambda(\beta)}.
\end{eqnarray}

\noindent Similarly, we write

\begin{eqnarray*}
 \overleftarrow{W}_{M,N}(x)=e^{-(N-M)\lambda}P^x [\exp \{ \beta \sum^{N-M}_{t=1} \eta(N-t,\omega_t)  \}].
\end{eqnarray*}

\noindent for the related backward expression. Another useful notation is the 'conditional' partition
function: take $M<N\in \ze$, $x,y \in \ze^d$,

\begin{eqnarray*}
 W_{M,N}(x|y)=e^{-(N-M)\lambda}P^x [\exp \{ \beta \sum^{N-M}_{t=1} \eta(M+t,\omega_t) \}|\, \omega_{N-M}=y].
\end{eqnarray*}

In Section $2$ we will prove Theorem \ref{convergence} and Theorem \ref{convergence2}, postponing the proof of Proposition \ref{concentration} to Section $3$.




\section{Proof of Theorems $2$ and $3$}
In this section, we make no specific assumptions on the environment. We just need Proposition \ref{concentration}
and Theorem \ref{invariance} to hold, which is the case for gaussian or bounded environments in the $L^2$ region.
Actually, we will see that we don't need the whole strength of the invariance principle, but just an averaged version
of it.

In both cases, the proof is a sequence of $L^2$ calculus. Let's begin describing the aforementioned Local
Limit Theorem as it appears in \cite{V}, page $6$, Theorem $2.3$: 
take $x,y\in \ze^d$, $M\in \ze$, then, for $A>0$, $N>M$ and $l_N=O(N^{\alpha})$ with $0<\alpha<1/2$,

\begin{eqnarray}\label{LLT}
 W_{M,N}(x|y)&=& {W}_{M,M+l_N}(x) \times \overleftarrow{W}_{N-l_N,N}(y)\times e^{\beta \eta(N,y)-\lambda(\beta)}\, \\
             \nonumber 
             &+& \, R_{M,N}(x,y),
\end{eqnarray}

\noindent where

\begin{eqnarray}\label{rest}
 \lim_{N\to +\infty} \sup_{|x-y|<A N^{1/2}} Q(R^2_{M,N}(x,y))=0.
\end{eqnarray}

\noindent Note that by symmetry, this result is still valid when we fix $y$ and take the
limit $M \to -\infty$. 

\vspace{3ex}
 
\noindent \textsc{ Proof of Theorem \ref{convergence}:} We first restrict the sum in (\ref{density}) to 
the region of validity of the local limit theorem. This is done using the quenched central
limit theorem averaged with respect to the disorder:

\begin{eqnarray}\label{separation}
 q_N = \sum_{|x|<A N^{1/2}} \mu^x_{-N,0}(\omega_N=0) + \sum_{|x|\geq AN^{1/2}}\mu^x_{-N,0}(\omega_N=0).
\end{eqnarray}

\noindent We compute the $L^1$ norm of the second term of the sum making use
of the invariance by translation under the law $Q$:

\begin{eqnarray}\label{queue_tcl}
 Q \sum_{|x|\geq AN^{1/2}}\mu^x_{-N,0}(\omega_N=0)= Q \mu^0_{0,N}(|\omega_N|\geq AN^{1/2}).
\end{eqnarray}

\noindent By Theorem \ref{invariance}, we have that  

\begin{eqnarray*}
\lim_{A\to +\infty} \limsup_{N\to +\infty} Q  \mu^x_{0,N}(|\omega_N| \geq A N^{1/2})=0.
\end{eqnarray*}

\smallskip

\noindent We can now concentrate on the first term in (\ref{separation}). We denote by $p(\cdot,\cdot)$ (resp. $p_N(\cdot,\cdot)$ )
the transition probabilities of the simple symmetric random walk on $\ze^d$  (resp.  
its $N$-step transition probabilities). Thanks to (\ref{LLT}),

\begin{eqnarray}\label{calc1}
 \nonumber
 \sum_{|x|<A N^{1/2}} \mu^x_{-N,0}(\omega_N=0)
\nonumber
 &=& \sum_{|x|<A N^{1/2}} \frac{W_{-N,0}(x|0)}{W_{-N,0}(x)}\,p_N(x,0)\\
 &=& \overleftarrow{W}_{-l_N,0}(0) e^{\beta \eta(0,0)-\lambda}
\!\!\! \sum_{|x|<A N^{1/2}} \frac{W_{-N,-N+l_N}(x)}{W_{-N,0}(x)}\,p_N(x,0)\\
\nonumber		
 &+&  \sum_{|x|<A N^{1/2}} \frac{R_{-N,0}(x,0)}{W_{-N,0}(x)}\,p_N(x,0).
\end{eqnarray}

 We again integrate the second term of the sum (\ref{calc1}). We use Cauchy-Schwarz inequality  and translation invariance:

\begin{eqnarray}\nonumber
 Q \sum_{|x|<A N^{1/2}} \frac{R_{-N,0}(x,0)}{W_{-N,0}(x)}\,p_N(x,0)
&\leq& \sup_{|x|<AN^{1/2}} \{ Q(R^2_{-N,0}(x,0))^{1/2} \} \\ \label{restineq}
&\times& \, Q(W^{-2}_{-N,0}(0))^{1/2} \, \sum_{|x|<A N^{1/2}} p_N(x,0).
\end{eqnarray}

\noindent The first term in the right side tends to zero thanks to (\ref{rest}),
the second one is easily seen to be bounded thanks to Proposition \ref{concentration} and
the third one is less than

\begin{eqnarray*}
\sum_x p_N(x,0)= \sum_x p_N(0,x)=1,
\end{eqnarray*}

\noindent so the left member of (\ref{restineq}) tends to zero. We are left to 
the study of the first summand in (\ref{calc1}),

\begin{eqnarray}\label{calc2}
 \nonumber
  &\ &  \overleftarrow{W}_{-l_N,0}(0)e^{\beta \eta(0,0)-\lambda}
\sum_{|x|<A N^{1/2}} \frac{W_{-N,-N+l_N}(x)}{W_{-N,0}(x)} \,p_N(x,0)\\
&\ & \quad = \quad \overleftarrow{W}_{-l_N,0}(0) e^{\beta \eta(0,0)-\lambda}\sum_{|x|<A N^{1/2}}p_N(x,0)\\
\nonumber
&\ &\quad \quad +\quad \overleftarrow{W}_{-l_N,0}(0) e^{\beta \eta(0,0)-\lambda}
\sum_{|x|<A N^{1/2}}\left\lbrace  \frac{W_{-N,-N+l_N}(x)}{W_{-N,0}(x)}-1 \right\rbrace p_N(x,0).
\end{eqnarray}

\noindent We see that we are done as long as we can control the convergence of the second summand in (\ref{calc2}).
It will be enough to prove that it converges to zero in probability. Let us denote

\begin{eqnarray*}
 g_N:= \overleftarrow{W}_{-l_N,0}(0) e^{\beta \eta(0,0)-\lambda},
\end{eqnarray*}

\begin{eqnarray*}
 h_N:= \sum_{|x|<A N^{1/2}}\left\lbrace  \frac{W_{-N,-N+l_N}(x)}{W_{-N,0}(x)}-1 \right\rbrace p_N(x,0).
\end{eqnarray*}

\noindent We already know that $\{ g_N:\, N\geq 1\}$ is bounded in $L^1(Q)$. It is then enough to
prove that $h_N$ tends to zero in $L^1$. Indeed, using translation invariance,
Cauchy-Schwarz inequality and the uniform boundedness of negative moments
of $W_{-N,0}$,

\begin{eqnarray*}
 Q(|h_N|) &\leq& \sum_{|x|<A N^{1/2}}Q \left( \left|  \frac{W_{-N,-N+l_N}(x)}{W_{-N,0}(x)}-1 \right| \right) p_N(x,0)\\
        &\leq& Q \left( \left|  \frac{W_{-N,-N+l_N}(0)}{W_{-N,0}(0)}-1 \right| \right)\\
        &\leq& Q \left( \left|W_{0,l_N}(0)-W_{0,N}(0)\right|^2 \right)^{1/2} Q\left( W^{-2}_{-N,0}(0)\right)^{1/2}.
\end{eqnarray*}

\noindent This clearly tends to zero. It is now a simple exercise to show that
$Q(g_N |h_N|\geq \delta)$ tends to zero as $N$ tends to infinity for all $\delta>0$.

\noindent So far, we have proved that $q_N$ tends to $q_{+\infty}$ in $Q$-probability.
But, by an elementary result, we know that, for $q_n,\, q_{+\infty}>0$, convergence in probability implies
$L^1$ convergence as long as $Q(|q_N|) \to Q(|q_{+\infty}|)$ (which is clearly the case here 
because all these expressions are equal to $1$).
This finishes the proof of Theorem \ref{convergence}. $\square$

\vspace{2ex}




\textsc{ Proof of Theorem \ref{convergence2}:} The details are very similar to
the previous proof. We split the sum in (\ref{densityNM}) according to $|x|\leq A N^{1/2}$ or not.
We can apply the central limit theorem to

\begin{eqnarray*}
Q \sum_{|x|>AN^{1/2}}\mu^x_{-N,M}(\omega_N=0)= Q \mu^0_{0,N+M}\left( 
|\omega_N|>AN^{1/2} \right).
\end{eqnarray*}

\noindent Now, by the Markov property and the local limit theorem, we have

\begin{eqnarray*}
\sum_{|x|\leq AN^{1/2}} \mu^x_{-N,M}\left(\omega_N=0 \right)&=&
\sum_{|x|\leq AN^{1/2}}\frac{W_{0,M}(0)\, W_{-N,0}(x|0)}{W_{-N,M}(x)}\, p_N(x,0)\\
&=& W_{0,M}(0) \overleftarrow{W}_{-l_N,0} 
    \sum_{|x|\leq AN^{1/2}} \frac{W_{-N,N+l_N}(x)}{W_{-N,M}(x)}p_N(x,0)\\
&+& W_{0,M}(0) \sum_{|x|\leq AN^{1/2}} \frac{R_{-N,0}(x,0)}{W_{-N,M}(x)}p_N(x,0)
\end{eqnarray*}

\noindent The second summand is again treated using the Cauchy-Schwarz inequality,
(\ref{rest}) and the independence of  $W_{0,M}(0)$ and $W^{-1}_{-N,M}(x)$. The
first summand as to be written as

\begin{eqnarray}.\nonumber
&\ &W_{0,M}(0) \overleftarrow{W}_{-l_N,0}(0) \sum_{|x|\leq AN^{1/2}} p_N(x,0)\\
\label{findec}
&+& W_{0,M}(0) \overleftarrow{W}_{-l_N,0}(0)
\sum_{|x|\leq AN^{1/2}} \left\lbrace \frac{W_{-N,N+l_N}(x)}{W_{-N,M}(x)}(0) - 1\right\rbrace p_N(x,0).
\end{eqnarray}

\noindent The second summand of (\ref{findec}) can be handled like the one in (\ref{calc2}), and using
the independence of $W_{0,M}(0)$ and  $\overleftarrow{W}_{-l_N,0}(0)$. $\square$

\vspace{3ex}




\section{Concentration inequalities} 
The proof follows closely \cite{T}, Section $2$ (see the proof of the lower bound of Theorem $1.1$ therein). 
Recall that we assumed that the environment is bounded by one.
In the $L^2$ region, it is known that $QZ^2_N \leq K (QZ_N)^2$ (see (\ref{L2})). This implies that

\begin{eqnarray}\label{PZ}
 Q \left( Z_N \geq \frac{1}{2}QZ_N \right) \geq \frac{1}{4} \frac{(QZ_N)^2}{QZ^2_N}\geq \frac{1}{4K},
\end{eqnarray}

\noindent thanks to Paley-Zigmund inequality.

The following is easily proved (\cite{CH}, proof of Theorem $1.5$): Let

$$A= \left\lbrace  Z_N \geq \frac{1}{2} QZ_N,\, \langle L_N(\omega,\tilde{\omega})\rangle^{(2)}_{N} \leq C \right\rbrace.$$

\noindent where the brackets mean expectation with respect to two independent copies of the polymer measure on the same
environment. Then, we can find $C>1$ such that, for all $N\geq 1$,
 
 \begin{eqnarray}\label{aff}
 Q(A)\geq 1/C.
 \end{eqnarray}

\noindent  It is convenient to see $Z_N$ as a function from $[-1,1]^{T_N}$ to $\re$, where 
$T_N=\{ (t,x):\, 0\leq t\leq N,\, |x|_1\leq t\}$ and $|x|_1=\sum^d_{i=1}|x_i|$ for $x=(x_1,...,x_d)$ . 
For $u>0$, let $B=\left\lbrace  z\in [-1,1]^{T_N}: \log Z_N(z)\leq \lambda N 
- \log 2 -u \right\rbrace $. This is a convex compact subset of $[-1,1]^{T_N}$.

\vspace{3ex}

In order to apply theorems for concentration of product measures, we need to introduce some
notation. For $x\in [-1,1]^{T_N}$, let

\begin{equation}
 U_B(x)= \{ h(x,y): y\in B \},
\end{equation}

\noindent where $h(x,y)_i={\bf 1}_{x_i \neq y_i}$. Let also $V_B(x)$ be the convex envelope $U_B(x)$ when
we look at it as a subset of $\re^{T_N}$. Finally, let $f(x,B)$ be the euclidean distance from the origin to
$V_B(x)$. Let us state the following result from \cite{Tnl}, Theorem $6.1$:

\begin{theorem}\label{Tal}
$$ \int \exp \left\lbrace  \frac{1}{4} f^2(x,B)\right\rbrace  \, dQ(x) \leq \frac{1}{Q(B)}.$$
\end{theorem}

\vspace{3ex}

On the other hand, for $x\in [-1,1]^{T_N},\, y\in B$, we have $|x_i-y_i|\leq 2h(x,y)_i$. Then for
every finite sequence $(y^{(k)})^M_{k=1} \subset B$ and any sequence of non-negative number such that
$\sum^M_{k=1} \alpha_k=1$, we have that 

$$|x_i - \sum^M_{k=1} \alpha_k y^{(k)}_i| \leq 2 \sum^M_{k=1}\alpha_k h(x,y^{(k)})_i,$$

\noindent for every $i \in T_N$. This yields

$$||x- \sum^M_{k=1} \alpha_k y^{(k)} || \leq 2 \left( \sum_{i\in T_N} \{ \sum^M_{k=1} \alpha_k
h(x,y^{(k)}_i) \}^2 \right)^{1/2}.$$

\noindent We can now optimize over all convex combinations of elements of $B$ (remember that it is convex),
we obtain

$$d(x,B) \leq 2 f(x,B),$$

\noindent where $d(x,B)$ is the euclidean distance from $x$ to $B$. We use Theorem \ref{Tal} to conclude
that

\begin{eqnarray*}
 \int_A \exp \left\lbrace  \frac{1}{16} d^2(x,B) \right\rbrace \, dQ(x) \leq \frac{1}{Q(B)}.
\end{eqnarray*}

\noindent We can find $\overline{x}\in A$ such that $Q(A) \exp \{ \frac{1}{16} d^2(\overline{x},B) \} \leq 1/Q(B)$. Using
(\ref{aff}), we find that $1/C \exp \{ \frac{1}{16} d^2(\overline{x},B) \} \leq 1/Q(B)$. Let $q=Q(B)$. By compacity
and some simple calculations, we conclude that we can find $\overline{z}\in B$ such that

\begin{eqnarray}\label{conc0}
 d(\overline{x},\overline{z}) \leq 4 \sqrt{\log \frac{C}{q}}.
\end{eqnarray}

\vspace{3ex}

Now, 

\begin{eqnarray*}
 Z_N(\overline{z})=Z_N(\overline{x}) \left\langle  \exp \left\lbrace 
                   \beta \sum^N_{t=1} (\overline{z}(t,\omega_t)-\overline{x}(t,\omega_t)) \right\rbrace \right\rangle _{\overline{x}}\\
                  \geq Z_N(\overline{x}) \exp \beta \left\langle  \sum^N_{t=1} (\overline{z}(t,\omega_t)-\overline{x}(t,\omega_t)) \right\rangle _{\overline{x}},
\end{eqnarray*}

\noindent where the brackets mean expectation with respect to the polymer measure in the $\overline{x}$ environment. But
using successively the Cauchy-Schwarz inequality, (\ref{conc0}) and the fact that $\overline{x}\in A$, we find that

\begin{eqnarray*}
 \left|\left\langle \sum^N_{t=1} (\overline{z}(t,\omega_t)-\overline{x}(t,\omega_t)\right\rangle_{\overline{x}} \right|
&=& \left|\sum^N_{t=1} \sum_a (\overline{z}(t,a)-\overline{x}(t,a) \left\langle  {\bf 1}_{\omega_t=a}\right\rangle_{\overline{x}}\right|\\
&\leq& ||\overline{x}-\overline{z}|| \left\langle L_N(\omega,\tilde{\omega})\right\rangle^{1/2}_{\overline{x}}\\
&\leq& 4 \sqrt{\log \frac{C}{q}} \sqrt{C}.
\end{eqnarray*}

\noindent So, using again the fact that $\overline{x}\in A$,

\begin{eqnarray*}
 \log Z_N(\overline{z}) &\geq& \log Z_N(\overline{x}) - 4 |\beta| \sqrt{\log \frac{C}{q}} \sqrt{C}\\
                        &\geq& \log \left( \frac{1}{2} QZ_N \right) - 4 |\beta| \sqrt{\log \frac{C}{q}} \sqrt{C}.
\end{eqnarray*}

\noindent Recalling that $\overline{z}\in B$ and after a few calculations, we conclude that

\begin{eqnarray*}
 q\leq C \exp \left\lbrace - \frac{u^2}{16C \beta^2}\right\rbrace, 
\end{eqnarray*}

\noindent which finishes the proof.

\smallskip
\section{}

\end{document}